\documentclass[a4paper,11pt]{article}
\usepackage{amsmath,amssymb,amsfonts,bm}
\usepackage[dvips]{graphicx}
\usepackage{theorem}
\makeatletter

\@addtoreset{equation}{section}
\makeatother
\newcommand{\R}{{\mathbb{R}}}
\renewcommand{\S}{{\mathbb{S}}}
\newcommand\QED
 {\hfill\hbox{\vrule width 4pt height 6pt depth 1.5pt}\par\bigskip}
\newcommand\argmax{\mathrm{argmax}}

\newcommand\dist{\mathrm{dist}}
\newcommand\proj{\mathrm{pr}}
\newcommand\sgn{\mathrm{sgn}}
\newcommand\spn{\mathrm{span}}

\newcommand\Tube{\mathrm{Tube}}
\newcommand\Vol{\mathrm{Vol}}

{\theorembodyfont{\rmfamily}\newtheorem{thm}{Theorem}[section]}
{\theorembodyfont{\rmfamily}\newtheorem{prop}{Proposition}[section]}
{\theorembodyfont{\rmfamily}\newtheorem{df}{Definition}[section]}
{\theorembodyfont{\rmfamily}\newtheorem{rem}{Remark}[section]}
\newenvironment{proof}[1][Proof]
 {\begin{trivlist}\item[\hskip\labelsep{\bfseries #1}]}{\end{trivlist}}
\title{
The tube method for the moment index in projection pursuit
}
\author{
Satoshi Kuriki\footnote{
Institute of Statistical Mathematics
and Graduate University for Advanced Studies, 
4-6-7 Minami-azabu, Minato-ku, Tokyo 106-8569, Japan, 
Email\,: {\tt kuriki@ism.ac.jp}
}\ \ %
and Akimichi Takemura\footnote{
Graduate School of Information Science and Technology,
University of Tokyo, 7-3-1 Hongo, Bunkyo-ku, Tokyo 113-0033, Japan,
Email\,: {\tt takemura@stat.t.u-tokyo.ac.jp}
}
}
\date{}
\begin{document}
\maketitle
\begin{abstract}
The projection pursuit index defined by a sum of squares of the third and
the fourth sample cumulants is known as the moment index proposed by
Jones and Sibson \cite{Jones-Sibson87}.
Limiting distribution of the maximum of the moment index under the null
hypothesis that the population is multivariate normal is shown to be
the maximum of a Gaussian random field with a finite Karhunen-Lo\`eve
expansion.
An approximate formula for tail probability of the maximum,
which corresponds to the $p$-value,
is given by virtue of the tube method through determining
Weyl's invariants of all degrees and the critical radius
of the index manifold of the Gaussian random field.


\end{abstract}

\section{Introduction}
\label{sec:introduction}

\subsection{Assessing the significance in projection pursuit}

Suppose that for each of $n$ individuals, a $q$ dimensional random
vector $x_t\in\R^q$, $t=1,\ldots,n$, is observed as an i.i.d.\ sample.
In the analysis of such multidimensional data, 
the projection of the $q$ dimensional data onto a lower dimensional
subspace is often used for the sake of interpreting the data.
In such cases it is important to select the subspace which clarifies
features of the data interesting to the analyst. 
In the principal components analysis or the canonical correlation
analysis, the subspaces are selected based on the variance of data %
 (\cite{Anderson03}).  The exploratory projection pursuit is the
method for detecting the subspace based on the non-normality of
data (\cite{Huber85}).
As a similar method, the Fast ICA (independent component analysis) %
is known (\cite{Hyvarinen-etal01}).

Let $\S^{q-1}$ be a set of $q$ dimensional unit vectors, or
the set of directional vectors in $\R^q$.
In the one dimensional projection pursuit, for each directional vector
 $h\in\S^{q-1}$, 
the projection pursuit index $I_n(h)$ is defined as a measure for the
non-normality of one dimensional projected data
\begin{equation}
\label{zt} 
 z_t=\langle x_t,h\rangle \in\R, \quad t=1,\ldots,n,
\end{equation}
and then 
the direction $h^*=\argmax\,I_n(h)$ attaining the maximum of the index
is searched.

However, the index $I_n(h)$ is a random function of $h$ depending on
the samples $x_t$'s.
Even when $x_t$'s were distributed according to the multidimensional
normal distribution, the function $I_n(h)$ is not constant, 
and the direction $h^*$ which achieves the maximum exists.
Therefore, it is important to assess whether it is not a pseudo peak
caused by stochastic fluctuations.  For this purpose, the framework of
the multiple testing can be employed.
Consider the null hypothesis that the data are distributed according to
the multidimensional normal distribution
\begin{equation}
\label{H0}
 H_0\,:\,x_t \sim N_q(\mu,\Sigma) \quad \rm i.i.d.,
\end{equation}
and let
\[
 \bar F_n(c) = P\Bigl(\max_{h\in \S^{q-1}}I_n(h) \ge c \mid H_0 \Bigr)
\]
be the upper probability of the maximum of $I_n(h)$ under the null hypothesis.
Then, $\bar F_n(I_n(h^*))$ is the $p$-value in the sense of
multiple testing, and we can use the $p$-value as a measure of
the significance of the maximum (Sun \cite{Sun91}).

Sun \cite{Sun91} described the limiting null distribution of
the maximum of Friedman \cite{Friedman87}'s projection pursuit index
in terms of a Gaussian random field as sample size goes to infinity,
and gave an approximation formula for it by an integral-geometric method
referred to as the tube method (\cite{Sun93}).
In this paper we gives an approximation formula for the moment index
proposed by Jones and Sibson \cite{Jones-Sibson87} by the tube method.

The moment index treated here is as follows:
Let $K_{k,n}(h)$ be the $k$th sample cumulant
of the projected data $z_t$ in (\ref{zt}),
and let 
$B_{1,n}(h)=K_{3,n}(h)/K_{2,n}(h)^{3/2}$ and  
$B_{2,n}(h)=K_{4,n}(h)/K_{2,n}(h)^2$
be the sample skewness and the sample kurtosis, respectively.
Then the moment index is defined by
\begin{equation}
\label{JS}
 I_n(h) = \frac{n}{6}B_{1,n}(h)^2 + \frac{n}{24} B_{2,n}(h)^2.
\end{equation}
Differently from Sun \cite{Sun91}'s treatment for Friedman's index,
we can determine geometric invariants of all degrees, 
and accordingly give an accurate formula for the $p$-values.

The structure of the paper is as follows:
The main results are summarized in Section \ref{sec:main}.
There, the limiting distribution of the maximum of
the moment index is described in terms of a Gaussian random field
with a finite Karhunen-Lo\`eve expansion, and determine
the geometric invariants of the index manifold.
An approximation formula for the upper probability of the maximum
can be obtained by incorporating these invariants.
Some numerical experiments to examine their accuracy are given there.
The main results of Section \ref{sec:main} are proved
in Section \ref{sec:proofs}.
Prior to Section \ref{sec:proofs}, we give a brief summary of the
tube method in Section \ref{sec:tube} as far as required.

\subsection{The tube method}

Here we give a very brief historical review of the tube method.

As explained in Section \ref{sec:tube}, the term tube
means a spherically tubular neighborhood around a set in the sphere.
Hotelling \cite{Hotelling39} pointed out a relation between the
$p$-value of a testing problem in nonlinear regression
and the volume of the tube, and demonstrated to calculate the $p$-value
by presenting a one dimensional formula for the volume of tube.
Weyl \cite{Weyl39} generalized Hotelling \cite{Hotelling39}'s formula
to the general dimensional case.  More recently,
Knowles and Siegmund \cite{Knowles-Siegmund89} and Sun \cite{Sun93} %
 found out the relation between the formula for the volume of tube
and the tail probability formula for the maximum of a Gaussian random field.
Since then, the tube method was applied to statistical
problems such as calculating null distributions of max-type test
statistics, or adjusting the multiplicity in multiple testing problems.
For example, the asymptotic distribution of the
Anderson-Stephens statistic (\cite{Anderson-Stephens72}) for
testing the uniformity of direction can be evaluated
 (\cite{Kuriki-Takemura04}).
For the other examples, see \cite{Kuriki-Takemura01} and \cite{Kuriki05}.
Nowadays, the tube method was proved 
by Takemura and Kuriki \cite{Takemura-Kuriki02} to be a special case of
the Euler characteristic heuristics,
which is known as another approach for approximating the distribution of
the maxima of random fields developed by 
Adler (\cite{Adler81}, \cite{Adler00}), 
Worsley (\cite{Worsley95as}, \cite{Worsley95aap})  
and Taylor  (\cite{taylor-hakuron}). 
For recent developments of the tube method and the Euler characteristic 
heuristic, see Adler and Taylor \cite{Adler-Taylor07}.
See also \cite{Kuriki-Takemura-ams}.

\section{Main results}
\label{sec:main}


We begin with giving the limiting distribution of the moment index
$I_n(h)$ in (\ref{JS}) under the null hypothesis of multivariate normality.
Without loss of generality, we assume that $x_t$'s are distributed
according to the $q$ dimensional standard normal distribution $N_q(0,I_q)$.
\begin{thm}
\label{thm:clt}
Let $\xi_1\in \R^{q^3}$, $\xi_2\in \R^{q^4}$ be random vectors consisting of
independent standard normal random variables.
For a unit vector $h\in \S^{q-1}$, let
\begin{equation}
\label{Z1Z2}
 Z_1(h) = \langle h\otimes h\otimes h,\xi_1\rangle, \quad
 Z_2(h) = \langle h\otimes h\otimes h\otimes h,\xi_2\rangle,
\end{equation}
where $\otimes$ denotes the Kronecker product.
Under the null hypothesis $H_0$ in (\ref{H0}), as $n\to\infty$, 
$\max_{h\in \S^{q-1}} I_n(h)$ converges in distribution to
$\max_{h\in \S^{q-1}} I(h)$, 
where
\[
 I(h) = Z_1(h)^2+Z_2(h)^2.
\]
\end{thm}
\begin{proof}
Let $C(\S^{q-1})$ be the Banach space of real valued continuous functions
on $\S^{q-1}$ endowed with the supremum norm. 
Note that the sample cumulant $K_{k,n}(\cdot)$, 
the sample skewness $B_{1,n}(\cdot)$,
the kurtosis $B_{2,n}(\cdot)$, and
the moment index $I_n(\cdot)$ are the elements of $C(\S^{q-1})$.
Theorem 2.1 of Kuriki and Takemura \cite{Kuriki-Takemura01} states that
the $(\sqrt{n} B_{1,n}(\cdot), \sqrt{n} B_{2,n}(\cdot))$ 
converges to $(\sqrt{3!} Z_1(\cdot), \sqrt{4!} Z_2(\cdot))$
in distribution in the space $C(\S^{q-1})$.  The theorem follows
from the continuous mapping theorem.
\QED
\end{proof}

By means of Theorem \ref{thm:clt} above, for large sample size $n$,
the $p$-value $\bar F_n(I_n(h^*))$ can be approximated by
$\bar F(I_n(h^*))$ with
\[
 \bar F(c) = P\Bigl(\max_{h\in \S^{q-1}} I(h) \ge c\Bigr).
\]
Moreover, by letting
\[
 Z(h,\theta) = \cos\theta Z_1(h) + \sin\theta Z_2(h), \quad
 h\in \S^{q-1},\ \theta \in\Bigl(-\frac{\pi}{2},\frac{\pi}{2}\Bigr]
\]
with $Z_1(h)$ and $Z_2(h)$ given in (\ref{Z1Z2}), we have
\begin{equation}
\label{max}
 \Bigl\{\max_{h\in \S^{q-1}} I(h)\Bigr\}^{1/2}
 = \max_{h\in \S^{q-1},\,\theta\in(-\pi/2,\pi/2]} Z(h,\theta).
\end{equation}
Therefore, from now on, we restrict our attention to the distribution
 of the maximum of $Z(h,\theta)$.


Let
\[
 p=q^3+q^4
\]
and
\begin{equation}
\label{M}
%
 M = \left\{
  (\cos\theta(h\otimes h\otimes h),\sin\theta(h\otimes h\otimes h\otimes h))
  \in \R^p \mid
  h\in \S^{q-1},\ \theta\in \Bigr(-\frac{\pi}{2},\frac{\pi}{2}\Bigl] \right\}.
\end{equation}
The maximum (\ref{max}) can be rewritten as
\begin{equation}
\label{Gaussian}
 \max_{x\in M} \langle x,\xi\rangle,
 \qquad \xi=(\xi_1,\ldots,\xi_p)\sim N_p(0,I_p).
\end{equation}
Note that $M$ and $\S^{q-1}\times (-\pi/2,\pi/2]$ are one-to-one.
It is easy to see that $M$ is a $q$ dimensional closed submanifold of
$\S^{p-1}$.
As shall be explained in Section \ref{sec:tube}, (\ref{Gaussian}) is
of the canonical form of the tube method in (\ref{Zx}).

The upper probability function of the chi-square distribution
with $\nu$ degrees of freedom is denoted by
\begin{equation}
\label{barG}
 \bar G_{\nu}(c)
 = \frac{1}{2^{\nu/2}\Gamma(\nu/2)} \int_c^\infty t^{\nu/2-1} e^{-t/2} dt.
\end{equation}
The volume of the $m-1$ dimensional volume of the unit sphere
$\S^{m-1}$ is denoted by
\begin{equation}
\label{Omega}
 \Omega_m=\frac{2\pi^{m/2}}{\Gamma(m/2)}.
\end{equation} 
The following is the main theorem of this paper.
The proof is given in Sections \ref{subsec:kappa} and \ref{subsec:rho_c}.
\begin{thm}
\label{thm:main}
As $c\to\infty$,
\begin{eqnarray*}
 P\Bigl(\max_{h\in \S^{q-1}} I(h) \ge c^2\Bigr)
= \sum_{e=0,\,e:{\rm even}}^q 
 \kappa_e \frac{\Gamma((q+1-e)/2)}{2^{1+e/2}\pi^{(q+1)/2}}
 \bar G_{q+1-e}(c^2)
 + O\bigl(c^{p-2}e^{-\rho_c c^2/2}\bigr),
\end{eqnarray*}
where 
\begin{equation}
\label{kappa}
 \kappa_e = \Omega_q\,
 \frac{(-3)^{e/2} (q-1)!}{(q-e)!} \sum_{j=0}^{e/2}
 \frac{(q-e-2j)}{(e/2-j)!\,j!} (-2)^j E_{(q-1-e)/2-j},
\end{equation}
\begin{equation}
\label{elliptic}
 E_k
= \int_{-\pi/2}^{\pi/2} (3 \cos^2\theta + 4\sin^2\theta)^k d\theta
\end{equation}
and
\begin{equation}
\label{rho_c}
 \rho_c = \frac{25}{16}.
\end{equation}
\end{thm}

\begin{rem}
$E_k$ in (\ref{elliptic}) with $k$ an integer or a half-integer
can be evaluated numerically by recurrence formulas: 
\begin{equation}
\label{forward}
 E_k = \frac{7(2k-1)}{2k} E_{k-1} - \frac{12(k-1)}{k} E_{k-2}, \quad
 \mbox{for\ \ } k=1,\frac{3}{2},2,\ldots,
\end{equation}
and
\begin{equation}
\label{backward}
 E_k = \frac{7(2k+3)}{24(k+1)} E_{k+1} - \frac{k+2}{12(k+1)} E_{k+2},
 \quad
 \mbox{for\ \ } k=-\frac{3}{2},-2,-\frac{5}{2},\ldots,
\end{equation}
with the boundary conditions
\[
E_{1/2} = 4 E(1/4) \doteq 4\times 1.46746, \quad
E_0 = \pi, \quad
E_{-1/2} = K(1/4) \doteq 1.68575, \quad
E_{-1} = \frac{\pi}{2\sqrt{3}},
\]
where
$E(1/4)$ and $K(1/4)$ are
complete elliptic integrals of the first kind and the second kind 
 (\cite{Abramowitz-Stegun92}, p.\,608--9).
The proofs for (\ref{forward}) and (\ref{backward}) are given in
 Section \ref{subsec:recurrence}.
\end{rem}


To conclude this section, we give numerical examples
for the purpose of examining the accuracy of the formula.
The tail probability of the maximum for $q=2$ is given by
\begin{equation}
\label{q2}
P\Bigl(\max_{h\in \S^{2-1}} I(h) \ge c^2\Bigr)
\,\sim\,w \bigl\{ \bar G_3(c^2) - \bar G_1(c^2) \bigr\}
\,=\,w \sqrt{\frac{2}{\pi}}\,c\,e^{-c^2/2}, \qquad c\to\infty,
\end{equation}
where $w = 2 E(1/4) \doteq 2\times 1.46746$.

\begin{figure}[h]
\begin{center}
\scalebox{0.5}{\includegraphics{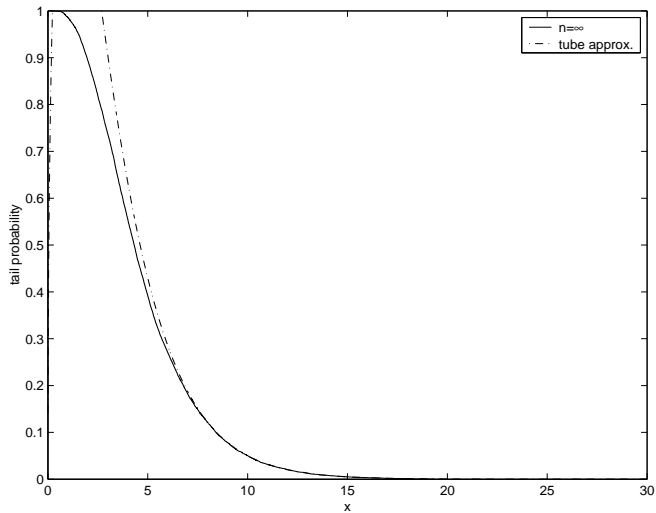}} 
\caption{Tail probability of limiting distribution (solid line)
 and its approximation by the tube method (dotted line).}
\label{fig:tube}
\end{center}
\end{figure}

Figure \ref{fig:tube} depicts the empirical upper probability of
the limiting distribution $P(\max_{h\in \S^{2-1}} I(h)\allowbreak \ge x)$
estimated by Monte Carlo simulations based on 10,000 replications,
and its approximation by the tube method.
One can see that the quantiles of the limiting distribution are
fully approximated by the tube method approximation (\ref{q2}).

\begin{figure}[h]
\begin{center}
\scalebox{0.5}{\includegraphics{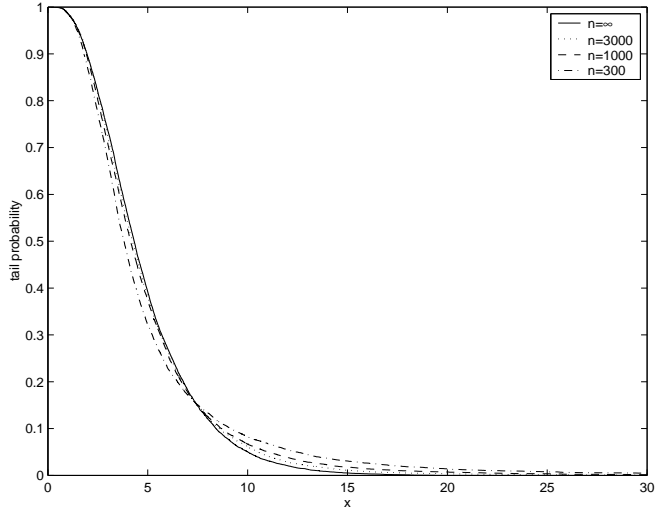}} 
\caption{Tail probabilities of finite sample distributions
 ($n=300,1000,3000,\infty$).}
\label{fig:finite}
\end{center}
\end{figure}

Figure \ref{fig:finite} depicts the empirical upper probability of
the finite sample distributions $P(\max_{h\in \S^{2-1}}\allowbreak I_n(h) \ge x)$
when $n=300,1000,3000,\infty$.  The number of replications is 10,000. 

\section{Summary of the tube method}
\label{sec:tube}

\subsection{Volume of the tubes and tail probabilities of the maxima}
\label{subsec:volume}


In this section we summarize the facts on the tube method
required for proving Theorem \ref{thm:main}.
We state Theorem \ref{thm:global} since its statement is 
not given in existing literature.

Let $\S^{p-1}$ be the unit sphere in $\R^p$, and let $M$ be a closed
 subset of $\S^{p-1}$.  Assume that $M$ is a $d$ dimensional
$C^2$ closed submanifold without boundaries embedded in $\S^{p-1}$,
and is endowed with the metric induced by the standard inner product
$\langle\cdot,\cdot\rangle$ of $\R^p$.

The set of points of $\S^{p-1}$ whose great circle distance (angle)
from $M$ is less than or equal to a constant $\theta$ is called 
the tube about $M$ with the radius $\theta$, and denoted by
\[
 \Tube(M,\theta) = \Big\{ y\in \S^{p-1} \mid \dist(y,M) \le \theta \Big\},
 \qquad \dist(y,M)=\min_{x\in M}\ \cos^{-1}\langle y,x\rangle.
\]
In a similar manner, the Euclidean tube is  defined
 in the Euclidean space by the usual distance.
But it does not play any role in this paper.

Let $y$ be a point of $\S^{p-1}\setminus M$.
The point $x=\proj(y)$ which attains the minimum $\min_{x\in M} \dist(y,x)$
is called the projection of $y$ onto $M$.
If $y$ is close to $M$, then $\proj(y)$ exists uniquely.
Whereas, if $y$ is far from $M$, then there can exist
two points $x_1,x_2\in M$ equidistant from $y$ which attain the minimum
$\min_{x\in M} \dist(y,x)$ simultaneously.
The supremum of the distances which assures the uniqueness
is called the critical radius.
\begin{df}
\label{df:critical}
When the $\proj(y)\in M$ is defined uniquely 
 for every $y \in \Tube(M,\theta)\setminus M$, 
it is said that the tube $\Tube(M,\theta)$ does not have a self-overlap.
The supremum
\[
\theta_c = \sup \{ \theta\ge 0 \mid
 \mbox{$\Tube(M,\theta)$ does not have a self-overlap} \} 
\]
is called the critical radius of $M$. 
\end{df}
The volume of a tube whose radius is less than or equal to the critical
radius $\theta_c$ can be calculated by taking a coordinate system
 based on the projection 
 (the Fermi coordinates).
The following proposition for the dimension $d=1$ is
due to Hotelling \cite{Hotelling39},
 and due to Weyl \cite{Weyl39} for the general dimensional case.
Here $\Omega_p$ denotes the $p-1$ dimensional volume of $\S^{p-1}$
defined in (\ref{Omega}),
and 
\[
 \bar B_{a,b}(c)
 = \frac{\Gamma(a+b)}{\Gamma(a)\Gamma(b)} \int_c^1 
t^{a-1} (1-t)^{b-1} dt
\]
is the upper probability function of the beta distribution
with parameters $(a,b)$.
\begin{prop}
\label{prop:weyl}
For $0\le\theta\le\theta_c$, $p-1$ dimensional volume of the tube is
 given by
\[
 \Vol(\Tube(M,\theta)) = \Omega_p \sum_{e=0,\ e:\,\rm even}^d \kappa_e
 J_e(\theta),
\]
where
\[
 J_e(\theta) = \frac{\Gamma((d+1-e)/2)}{2^{1+e/2}\pi^{(d+1)/2}}
 \bar B_{(d+1-e)/2,(p-d-1+e)/2}(\cos^2\theta),
\]
and the $\kappa_e$ is the intrinsic invariant of the manifold $M$
defined below in (\ref{invariant}), referred to as Weyl's curvature invariant.
\end{prop}

Let $\xi=(\xi_1,\ldots,\xi_p)$ be a random vector consisting of
independent standard normal random variables.
That is, $\xi\sim N_p(0,I_p)$.
Define a Gaussian random field on a submanifold $M$ of $\S^{p-1}$ by
\begin{equation}
\label{Zx}
Z(x)=\langle x,\xi\rangle, \quad x\in M \ \ (\subset \S^{p-1}).
\end{equation}
This is a canonical form of Gaussian random fields
of mean 0 and variance 1 with a finite Karhunen-Lo\`eve expansion.

By replacing $\Omega_p\bar B$ by
the upper probability function of the chi-square distribution
$\bar G$ in (\ref{barG}), we have
an approximation formula for the tail probability
of the maximum of $Z(x)$ (\cite{Kuriki-Takemura01}, \cite{Takemura-Kuriki02}).
\begin{prop}
\label{prop:tail}
As $c\to\infty$,
\[
 P\Bigl(\max_{x\in M} Z(x)\ge c\Bigr)
 = \sum_{e=0,\,e:\rm even}^{d} \kappa_e \psi_e(c)
 + O\bigl(c^{p-2} e^{-(1+\tan^2\theta_c)c^2/2}\bigr),
\]
where
\[
 \psi_e(c) = \frac{\Gamma((d+1-e)/2)}{2^{1+e/2}\pi^{(d+1)/2}}
 \bar G_{d+1-e}(c^2).
\]
\end{prop}

Note that the larger the critical radius $\theta_c$ is, the smaller
the order of the remainder term is.

\subsection{Weyl's curvature invariants}
\label{subsec:weyl}

As we saw in Propositions \ref{prop:weyl} and \ref{prop:tail},
the Weyl's curvature invariants $\kappa_e$ and
the critical radius $\theta_c$ of the manifold $M$ are needed
in applying the tube method.
We will explain the way to determine them in this and subsequent subsections.

Write a local coordinate system of a $d$ dimensional closed manifold
$M$ as $(t^i)$.  The metric tensor is denoted by $g_{ij}$, 
and write the $(i,j)$th elements of the inverse of the
$d\times d$ matrix $(g_{ij})$
as $g^{ij}$.  Abbreviate $\partial/\partial t^i$ to $\partial_i$.
The connection coefficients and the curvature tensor are
given by
\begin{equation}
\label{Gamma}
 \Gamma_{ij}^k = \sum_{h=1}^d g^{kh} \Gamma_{ij,h}\quad
 \mbox{with}\quad
 \Gamma_{ij,k} = \frac{1}{2} (
   \partial_i g_{jk}
 + \partial_j g_{ik}
 - \partial_k g_{ij}),
\end{equation}
and
\begin{equation}
\label{R}
 R^{kl}_{ij}=\sum_h g^{lh}{R^k}_{hij}\quad
 \mbox{with}\quad
 {R^l}_{kij} = \partial_i\Gamma^l_{jk} - \partial_j\Gamma^l_{ik}
 + \sum_{s=1}^d (\Gamma^l_{is}\Gamma^s_{jk} - \Gamma^l_{js}\Gamma^s_{ik}),
\end{equation}
respectively.
Let
\begin{equation}
\label{H}
H^{kl}_{ij} = R^{kl}_{ij} - (\delta^k_i\delta^l_j - \delta^l_i\delta^k_j),
\end{equation}
where $\delta_i^j$ is Kronecker's delta.
For
$e=0,2,\ldots,[d/2]\times 2$, let
\begin{equation}
\label{He}
H_e = \sum_i \sum_\sigma \sgn(\sigma)\,
 H_{i_1 i_2}^{i_{\sigma(1)} i_{\sigma(2)}}
 H_{i_3 i_4}^{i_{\sigma(3)} i_{\sigma(4)}}
 \cdots
 H_{i_{e-1} i_e}^{i_{\sigma(e-1)} i_{\sigma(e)}}.
\end{equation}
Here the summation $\sum_i$ is taken over all sets
of $e/2$ paring made of distinct elements of $\{1,2,\ldots,d\}$,
that is, all possible ways of
$\{i_1,i_2,\ldots,i_e\}\subset \{1,2,\ldots,d\}$ satisfying
$i_1<i_2$, $i_3<i_4,\ldots,i_{e-1}<i_e$ and $i_1<i_3<\cdots<i_{e-1}$.
The summation $\sum_\sigma$ is 
taken over all permutations $\sigma$ of $\{1,2,\ldots,e\}$ such that
$\sigma(1)<\sigma(2)$, $\sigma(3)<\sigma(4),\ldots,\sigma(e-1)<\sigma(e)$.
Then, Weyl's curvature invariants are defined by
\begin{equation}
\label{invariant}
 \kappa_e = \int_M H_e\,\det(g_{ij})^{1/2} dt^1\cdots dt^d,
 \qquad e=0,2,\ldots,[d/2]\times d
\end{equation}
 (Weyl \cite{Weyl39}).

For instance, $H_e$ for $e=0,2,4$ are given as follows:
$H_0=1$, and hence $\kappa_0$ is the $d$ dimensional volume of $M$.
\[
 H_2 = \sum_{1\le i<j\le d} H^{ij}_{ij}
       = \frac{1}{2} \sum_{i,j=1}^d
 H^{ij}_{ij}
       = \frac{1}{2}\Bigl\{ \sum_{i,j=1}^d R^{ij}_{ij} - d(d-1) \Bigr\},
\]
where $\sum_{i,j=1}^d R^{ij}_{ij}$ is the scalar curvature.
\begin{eqnarray*}
H_4 &=& \sum_{1\le i<j<k<l\le d}
  (H_{ij}^{ij}H_{kl}^{kl} -H_{ij}^{ik}H_{kl}^{jl}
  +H_{ij}^{il}H_{kl}^{jk} +H_{ij}^{jk}H_{kl}^{il}
  -H_{ij}^{jl}H_{kl}^{ik} +H_{ij}^{kl}H_{kl}^{ij} \\
&& \qquad\qquad\quad
  -H_{ik}^{ij}H_{jl}^{kl} +H_{ik}^{ik}H_{jl}^{jl}
  -H_{ik}^{il}H_{jl}^{jk} -H_{ik}^{jk}H_{jl}^{il}
  +H_{ik}^{jl}H_{jl}^{ik} -H_{ik}^{kl}H_{jl}^{ij} \\
&& \qquad\qquad\quad
  +H_{il}^{ij}H_{jk}^{kl} -H_{il}^{ik}H_{jk}^{jl}
  +H_{il}^{il}H_{jk}^{jk} +H_{il}^{jk}H_{jk}^{il}
  -H_{il}^{jl}H_{jk}^{ik} +H_{il}^{kl}H_{jk}^{ij} ) \\
&=& \frac{1}{8}\sum_{i,j,k,l=1}^d
  (H_{ij}^{ij}H_{kl}^{kl} - 4 H_{ij}^{il} H_{kl}^{kj} + H_{ij}^{kl}H_{kl}^{ij})
%
\\
&=& \frac{1}{8}\Bigl\{
 \Bigl(\sum_{i,j=1}^d R_{ij}^{ij}\Bigr)^2
-4\sum_{i,j,k,l=1}^d R_{ij}^{il}R_{kl}^{kj}
 +\sum_{i,j,k,l=1}^d R_{ij}^{kl}R_{kl}^{ij} \\
&& \qquad\quad
 -2 (d-2)(d-3)\sum_{i,j=1}^d R_{ij}^{ij} + d(d-1)(d-2)(d-3) \Bigr\}.
\end{eqnarray*}

See Gray (\cite{Gray04}, Lemma 4.2) for the invariants of a Euclidean tube.

\subsection{Evaluation of critical radius}
\label{subsec:critical}

In this subsection we give theorems useful in calculating 
the critical radius of a closed submanifold of the sphere.
\begin{prop}
\label{prop:critical}
The critical radius $\theta_c$ of a closed submanifold $M$
of $\S^{p-1}$ satisfies
\begin{equation}
\label{critical}
 \cot^2\theta_c
 = \sup_{y,x\in M,\,y\ne x}
 h(x,y), \qquad h(x,y)=
   \frac{1-\langle y,P_x y\rangle}{(1-\langle x,y\rangle)^2},
\end{equation}
where $P_x$ is the orthogonal projection onto
the linear subspace $\spn\{x\}\oplus T_x M$ of $\R^p$, 
and $T_x M$ is the tangent space of $M$ at $x$
 (\cite{Johansen-Johnstone90}, \cite{Kuriki-Takemura01}).
\end{prop}

A theorem corresponding to a Euclidean tube is given by
Federer (\cite{Federer59}, Theorem 4.18).

The radius $\theta_c^{\rm loc}$ satisfying
\begin{equation}
\label{local}
 \cot^2\theta_c^{\rm loc}
 = \limsup_{y,x\in M,\,\Vert y-x\Vert\to 0}
 h(x,y)
\end{equation}
is called the local critical radius,
which is characterized as the curvature radius of $M$ at $x$ %
 (\cite{Johansen-Johnstone90}, \cite{Kuriki-Takemura01}).
By definitions, $\theta_c^{\rm loc}\ge \theta_c$,
and the equality holds if the supremum in (\ref{critical}) is attained
when $\Vert y-x\Vert\to 0$.

Define a real-valued function on $M\times M$ by
$r(x,y)=\langle x,y\rangle$.
This is the covariance function of the Gaussian random field (\ref{Zx}).
Denote the local coordinate system about $x$ and $y$
by $(s^i)$, $(t^i)$, respectively.

The set of the critical points of $r(x,y)$ which are not contained
in the diagonal set is denoted by
\[
 C = \Bigl\{ (x,y) \in M\times M \mid x\ne y,\ %
 \frac{\partial}{\partial s^i} r(x,y) = 0,\ %
 \frac{\partial}{\partial t^i} r(x,y) = 0 \Bigr\}.
\]
Then we have the following theorem.
\begin{thm}
\label{thm:global}
The critical radius $\theta_c$ satisfies
\[
 \theta_c = \min \Bigl\{ \theta_c^{\rm loc}, 
\inf_{(x,y)\in C} \frac{1}{2}\cos^{-1} \langle x,y\rangle \Bigr\}.
\]
\end{thm}

\begin{proof}
By Lemma 5.2 of \cite{Taylor-etal05}, 
%
if the supremum of $h(x,y)$
is attained at a point not contained in the diagonal set, then it
 belongs to $C$.  Furthermore, for the points $(x,y)\in C$,
it holds that $P_x y=\langle x,y\rangle x$,
\[
 h(x,y)
 = \frac{1-\langle x,y\rangle^2}{(1-\langle x,y\rangle)^2}
 = \frac{1+\langle x,y\rangle}{1-\langle x,y\rangle}
 = \cot^2 \Bigl( \frac{1}{2}\cos^{-1}\langle x,y\rangle \Bigr),
\]
and hence
\[
  \sup_{(x,y)\in C} h(x,y)
= \cot^2 \Bigl( \inf_{(x,y)\in C}\frac{1}{2}\cos^{-1}\langle x,y\rangle \Bigr).
\]
Since the supremum of $h(x,y)$ over the diagonal set is
 $\cot^2\theta_c^{\rm loc}$,
the theorem follows from Proposition \ref{prop:critical}.
\QED
\end{proof}

A theorem corresponding to a Euclidean tube with the dimension $d=1$
 is given by
Johansen and Johnstone (\cite{Johansen-Johnstone90}, Proposition 4.2).

\section{Proof of Theorem \ref{thm:main}}
\label{sec:proofs}

\subsection{Proof of (\ref{kappa})}
\label{subsec:kappa}


In this section, we prove Theorem \ref{thm:main}.
By means of Proposition \ref{prop:tail}, the approximation formula
for the upper probability of the maximum can be given through
determining Weyl's curvature invariants $\kappa_e$ and 
the critical radius $\theta_c$ of the index manifold $M$ in (\ref{M}).
The former is given here, and the latter is given in the next subsection.

The metric tensor, the connection coefficients, and the curvature tensor
 for $M$ are denoted by $g$, $\Gamma$, and $R$, respectively, as in
 Section \ref{subsec:weyl}.
Also, the same quantities for $\S^{q-1}$ are denoted by
$\bar g$, $\bar\Gamma$, and $\bar R$, respectively.

Write an element $h$ of $\S^{q-1}$ by a local coordinate system
as $h=h(t)$, $t=(t^i)$.
Let $h_i=\partial h/\partial t^i$.
The metric of $\S^{q-1}$ is $\bar g_{ij} = \langle h_i,h_j\rangle$.

An element $x$ of $M$ can be written as
\[
 x=
 (\cos\theta (h\otimes h\otimes h), \sin\theta (h\otimes h\otimes h\otimes h))
 \in M
\]
in terms of $(t,\theta)$.  The bases of the tangent space of $M$ are
\begin{eqnarray*}
&& \frac{\partial x}{\partial t^i} = (
      \cos\theta (h_i\otimes h\otimes h
                 +h\otimes h_i\otimes h
                 +h\otimes h\otimes h_i
                   ), \\
&& \qquad\quad\ %
       \sin\theta (h_i\otimes h\otimes h\otimes h
                  +h\otimes h_i\otimes h\otimes h
                  +h\otimes h\otimes h_i\otimes h \\
&& \qquad\qquad\qquad\qquad
                  +h\otimes h\otimes h\otimes h_i
                   )
     ), \quad
 i=1,\ldots,q-1, \\
&&
 \frac{\partial x}{\partial\theta} = (
      -\sin\theta (h\otimes h\otimes h),
       \cos\theta (h\otimes h\otimes h\otimes h\otimes h)
     ).
\end{eqnarray*}
In the following, $\theta$ is regarded as the 0th coordinate $t^0$ of $t$.
The metric tensor of $M$ is
\[
g_{ij} = \begin{cases}
 v(\theta) \bar g_{ij}(t) & \mbox{if $i,j\ne 0$}, \\
                        1 & \mbox{if $i=j=0$}, \\
                        0 & \mbox{otherwise}, \end{cases}
\]
where
\begin{equation}
\label{v}
 v(\theta) = 3\cos^2\theta+4\sin^2\theta
  = 3+\sin^2\theta = 4-\cos^2\theta.
\end{equation}
%
%
From this, the volume element of $M$ is shown to be
\begin{equation}
\label{dvol}
 \det(\bar g_{ij}(t))^{1/2} dt^1\cdots dt^{q-1}\,v(\theta)^{(q-1)/2}\,d\theta.
\end{equation}
Note that $\det(\bar g_{ij}(t))^{1/2} dt^1\cdots dt^{q-1}$
 is the volume element of $\S^{q-1}$.

Let $\dot v$ and $\ddot v$ be the first and second derivatives of
$v=v(\theta)$.
After some calculations along the lines with (\ref{Gamma}),
it is shown that
the non-zero connection coefficients of $M$ are
\[
 \Gamma_{ij}^k = \bar \Gamma_{ij}^k, \qquad
 \Gamma_{i0}^k =  \Gamma_{0i}^k =
   \frac{1}{2} \frac{\dot v}{v} \delta_i^k, \qquad
 \Gamma_{ij}^0 = - \frac{1}{2} \dot v \bar g_{ij} \qquad
 (i,j,k\ne 0),
\]
and all of the other coefficients are 0.

Next we will derive the curvature tensor by (\ref{R}).
Put
\[
 J^{kl}_{ij} = \delta^k_i\delta^l_j - \delta^l_i\delta^k_j.
\]
Noting that the curvature tensor of the unit sphere $\S^{q-1}$
is $\bar R^{kl}_{ij} = J^{kl}_{ij}$,  after cumbersome calculations
we see that the non-zero elements are
\[
 R^{kl}_{ij} = \Bigl\{ \frac{1}{v}
                -\frac{1}{4}\Bigl(\frac{\dot v}{v}\Bigr)^2 \Bigr\}
   J^{kl}_{ij}, \qquad
 R^{k0}_{i0} = -R^{0k}_{i0} = -R^{k0}_{0i} = R^{0k}_{0i} =
 \Bigl\{-\frac{1}{2}\frac{\ddot v}{v}
  +\frac{1}{4}\Bigl(\frac{\dot v}{v}\Bigr)^2 \Bigr\}
  \delta^k_i \qquad (i,j,k,l\ne 0).
\]
Furthermore, noting that
$\dot{v}=2\cos\theta\sin\theta$, 
$(\dot{v})^2=4\cos^2\theta\sin^2\theta=4(4-v)(v-3)=-4(v^2-7v+12)$,
$\ddot{v}=2\cos^2\theta-2\sin^2\theta=2(4-v)-2(v-3)=-2(2v-7)$,
%
%
we have the non-zero elements of $H^{kl}_{ij}$ in (\ref{H}) as
\[
 H^{kl}_{ij} = \alpha J^{kl}_{ij}, \qquad
 H^{k0}_{i0} = -H^{0k}_{i0} = -H^{k0}_{0i} = H^{0k}_{0i} =
  \beta \delta^k_i \qquad (i,j,k,l\ne 0),
\]
where
\[
 \alpha=\alpha(\theta) = -\frac{6}{v} + \frac{12}{v^2}, \qquad
 \beta=\beta(\theta) = -\frac{12}{v^2}.
\]
We substitute these quantities into (\ref{He}) to obtain
$H_e$, $e=0,2,\ldots,[q/2]\times 2$.

(i)\ \ The case where the set of the indices $i_1,i_2,\ldots,i_e$
 in the right-hand side of (\ref{He}) does not contain 0.
Because the number of the ways to make $e/2$ pairs from $q-1$ distinct
 objects is 
\[
 \frac{(q-1)!}{(q-1-e)!2^{e/2}(e/2)!},
\]
the summation of all terms corresponding to the case (i) becomes
\begin{equation}
\label{case-i}
 \alpha^{e/2}\times \frac{(q-1)!}{(q-1-e)!2^{e/2}(e/2)!}.
\end{equation}

(ii)\ \ The case where the set of the indices $i_1,i_2,\ldots,i_e$
 in the right-hand side of (\ref{He}) contains 0.
In this case, $i_1=0$, and
 $i_2\ne 0$, $\sigma(1)=1$ ($i_{\sigma(1)}=0$).
Noting that there are $q-1$ ways for $i_2$ ($i_2=1,\ldots,q-1$), and
that  $i_3,i_4,\ldots,i_e$ are indices resulting from making  
$e/2-1$ pairs from the set $\{1,2,\dots,q-1\}\setminus\{i_2\}$
having $q-2$ elements, 
the summation of all terms corresponding to the case (ii) becomes
\begin{equation}
\label{case-ii}
  (q-1)\beta\alpha^{e/2-1}\times
 \frac{(q-2)!}{(q-e)!2^{e/2-1}(e/2-1)!}.
\end{equation}
Summing up (\ref{case-i}) and (\ref{case-ii}) along with
\[
 \alpha^{e/2} = \Bigl(-\frac{6}{v}\Bigr)^{e/2}
  \sum_{j=0}^{e/2}{e/2\choose j}\Bigl(-\frac{2}{v}\Bigr)^j
\]
and
\[
 \beta\alpha^{e/2-1}
 = -\Bigl(-\frac{6}{v}\Bigr)^{e/2}
 \sum_{j=1}^{e/2}{e/2-1\choose j-1}\Bigl(-\frac{2}{v}\Bigr)^j
\]
yields
\[
 H_e = \sum_{j=0}^{e/2} \frac{A_j}{v^{e/2+j}},
\]
where
\[
 A_0 = \frac{(-3)^{e/2} (q-1)!}{(q-e-1)!(e/2)!},
\]
and for $j\ne 0$,
\begin{eqnarray*}
A_j
&=& (-6)^{e/2} (-2)^j \biggl\{
 \frac{(q-1)!}{(q-1-e)!2^{e/2}(e/2)!}{e/2 \choose j} 
 -\frac{(q-1)(q-2)!}{(q-e)!2^{e/2-1}(e/2-1)!}{e/2-1 \choose j-1}
\biggr\} \\
&=&  \frac{(-3)^{e/2} (q-1)!}{(q-e)!}
\frac{(q-e+2j) (-2)^j}{(e/2-j)!j!}.
\end{eqnarray*}

Since the expression for $A_j$ with $j\ne 0$ is consistent with
that for $A_0$, we have
\[
 H_e = \frac{(-3)^{e/2} (q-1)!}{(q-e)!} \sum_{j=0}^{e/2}
\frac{(q-e-2j)(-2)^j}{(e/2-j)!j!}\frac{1}{v^{e/2+j}}.
\]
Finally we obtain  $\kappa_e$ in (\ref{kappa})
by integrating $H_e$ over $M$
with respect to the volume element (\ref{dvol}).


\subsection{Proof of (\ref{rho_c})}
\label{subsec:rho_c}


\newcommand\tx{{\tilde x}}
\newcommand\ttheta{{\tilde\theta}}
\renewcommand\th{{\tilde h}}

In this subsection, making use of Theorem \ref{thm:global},
we show that the critical radius of the index manifold $M$ in (\ref{M}) is
$\theta_c = \tan^{-1}(3/4)$.
This implies that $\rho_c=1+\tan^2\theta_c=25/16$.
Throughout this subsection, we assume that vectors are column vectors
for notational convenience.
For instance, $\langle x,y\rangle=x'y$, where $'$ denotes the transpose.

We begin with obtaining the local critical radius $\theta_c^{\rm loc}$
by (\ref{local}).
Let
\[
 x=
 \begin{pmatrix}
   \cos\theta (h\otimes h\otimes h) \\
   \sin\theta (h\otimes h\otimes h\otimes h)
 \end{pmatrix}, \qquad
 \tx=
 \begin{pmatrix}
   \cos\ttheta (\th\otimes \th\otimes \th) \\
   \sin\ttheta (\th\otimes \th\otimes \th\otimes \th)
 \end{pmatrix}
\]
be two points of $M$.
Write for simplicity
$h_i=\partial h/\partial t^i$, 
$x_i=\partial x/\partial t^i$, $x_0=\partial x/\partial\theta$,
and $\bar G=(\bar g_{ij})$, $v=v(\theta)$ defined in (\ref{v}).
The orthogonal projection matrix onto $\spn\{x\}\oplus T_x M$ is
denoted by $P_x$.  Since $\spn\{x\}$ is orthogonal to $T_x M$, we have
\begin{eqnarray*}
\tx'P_x \tx &=&
 (\tx'x,\tx'x_1,\ldots,\tx'x_{q-1},\tx'x_0)
 \begin{pmatrix}1 & & \\ & v\bar G & \cr & & 1\end{pmatrix}^{-1}
 \begin{pmatrix}\tx'x \\ \tx'x_1 \\ \vdots \\ \tx'x_{q-1}
             \\ \tx'x_0 \end{pmatrix} \\
 &=& (\tx'x)^2
  + (\tx'x_1,\ldots,\tx'x_{q-1}) (v\bar G)^{-1}
    \begin{pmatrix}\tx'x_1 \cr \vdots \\ \tx'x_{q-1}\end{pmatrix}
  + (\tx'x_0)^2.
\end{eqnarray*}
The first term of the right-hand side is the square of
\[
\tx'x
= (\th'h)^3 \cos\ttheta\cos\theta
   +(\th'h)^4 \sin\ttheta\sin\theta.
\]
Noting that
\[
\tx'x_i
= w \th'h_i, \qquad
 w = 3(\th'h)^2 \cos\ttheta\cos\theta
    +4(\th'h)^3 \sin\ttheta\sin\theta,
\]
the second term becomes
\[
 \frac{w^2}{v}\,\th'(h_1,\ldots,h_{q-1})\bar G^{-1}
 \begin{pmatrix}h_1' \\ \vdots \\ h_{q-1}'\end{pmatrix} \th'
 = \frac{w^2}{v}\, \th'(I_q - hh') \th = \frac{w^2}{v}\,(1-(\th'h)^2).
\]
The third term is the square of
\[
\tx'x_0
= -(\th'h)^3 \cos\ttheta\sin\theta
    +(\th'h)^4 \sin\ttheta\cos\theta.
\]
Summing up these three terms,
the numerator of the right-hand side of (\ref{local}) is
\begin{eqnarray*}
1-\tx'P_x\tx
&=& 1- ((\th'h)^3 \cos\ttheta\cos\theta+(\th'h)^4 \sin\ttheta\sin\theta)^2
 - \frac{w^2}{v}\, (1-(\th'h)^2)  \\
&&  - (-(\th'h)^3 \cos\ttheta\sin\theta +(\th'h)^4 \sin\ttheta\cos\theta)^2 \\
&=& 1 -\cos^6\psi\cos^2\tilde\theta-\cos^8\psi\sin^2\tilde\theta \\
&&  - \frac{(3 \cos^2\psi \cos\tilde\theta\cos\theta
  +4 \cos^3\psi \sin\tilde\theta\sin\theta)^2}
           {3\cos^2\theta+4\sin^2\theta}\sin^2\psi \\
&=& f \ \ \mbox{(say)},
\end{eqnarray*}
where $\th'h=\cos\psi$.
On the other hand, the denominator of the right-hand side of (\ref{local}) is
\begin{eqnarray*}
 (1-\tx'x)^2
&=&
 (1-\cos^3\psi \cos\tilde\theta\cos\theta-\cos^4\psi \sin\tilde\theta\sin\theta)^2 \\
&=& g \ \ \mbox{(say)}.
\end{eqnarray*}


The local critical radius $\theta_c^{\rm loc}$ can be obtained by 
$\cot^2\theta_c^{\rm loc}=\limsup f/g$
when $\tilde\theta-\theta\to 0$, $\psi\to 0$.
Let $\tilde\theta-\theta=\delta$ and $u=\sin^2\theta$.
Ignoring $\psi^4$, $\psi^2\delta^2$, and $\delta^4$ as infinitesimals,
we have with aid of symbolic calculation that
\[
 f \sim 3(1+u)\psi^4 + \frac{12}{3+u}\psi^2\delta^2
\quad\mbox{and}\quad
 g \sim \frac{(3+u)^2}{4}\psi^4 + \frac{3+u}{2}\psi^2\delta^2
 + \frac{1}{4}\delta^4.
\]
Letting $\delta^2\sim k\psi^2$ for a constant $k$ (may be 0 or $\infty$),
we have
\[
 \frac{f}{g}\sim 12 \frac{(1+u)(3+u)+4k}{(3+u)(k+3+u)^2}.
\]
As a function of $k$, the right-hand side of the above takes its maximum 
\[
 \frac{48}{(3+u)^2(3-u)}
\]
at $k=(3+u)(1-u)/2$.
Furthermore as a function of $u$, this takes its maximum $48/27=16/9$ at $u=0$ over $0\le u\le 1$.
Note that when $u=0$, $k=3/2$ and $\theta=0$.

Summarizing the above arguments, one can see that
$\limsup f/g=16/9$ 
is attained when $\tilde\theta,\theta\to 0$, $\psi\to 0$,
$|\tilde\theta-\theta| \sim \sqrt{3/2} \psi$,
and accordingly
\[
 \theta_c^{\rm loc}
 = \cot^{-1} (\sqrt{16/9}) =\tan^{-1} (3/4) 
 \doteq 0.205 \pi.
\]



As the second step, we confirm that the local critical radius
is really the critical radius.
The covariance function of (\ref{Gaussian}) is
\begin{eqnarray*}
x'\tx
&=& \cos\theta \cos\ttheta (h'\th)^3
  + \sin\theta \sin\ttheta (h'\th)^4 \\
&=& \cos\theta \cos\ttheta \cos^3\psi
  + \sin\theta \sin\ttheta \cos^4\psi \\
&=& r(\psi,\theta,\ttheta) \ \ \mbox{(say)}.
\end{eqnarray*}
The ranges of the variables are
\begin{equation}
\label{range}
 \psi\in [0,\pi], \quad \theta,\ttheta \in
  \Bigl(-\frac{\pi}{2},\frac{\pi}{2}\Bigr].
\end{equation}
The set of the critical points are the set of the solutions of
\begin{eqnarray}
\label{eq1}
0 &=& \frac{\partial r}{\partial\psi}
= -\sin\psi ( 3 \cos\theta \cos\ttheta \cos^2\psi
  + 4 \sin\theta \sin\ttheta \cos^3\psi ), \\
\label{eq2}
0 &=& \frac{\partial r}{\partial\theta}
= -\sin\theta \cos\ttheta \cos^3\psi
  + \cos\theta \sin\ttheta \cos^4\psi, \\
\label{eq3}
0 &=& \frac{\partial r}{\partial\ttheta}
= -\cos\theta \sin\ttheta \cos^3\psi
  + \sin\theta \cos\ttheta \cos^4\psi.
\end{eqnarray}

(i)\ \ The case $\sin\psi\ne 0$.
From (\ref{eq1}), $(3 \cos\ttheta \cos^2\psi, 4 \sin\ttheta \cos^3\psi)$
is orthogonal to $(\cos\theta,\sin\theta)$.
From (\ref{eq2}), $(\cos\ttheta \cos^3\psi, \sin\ttheta \cos^4\psi)$
is orthogonal to $(-\sin\theta,\cos\theta)$.  Combining these,
\[
0
= 3 \cos^2\ttheta \cos^5\psi + 4 \sin^2\ttheta \cos^7\psi
= \cos^5\psi (3 \cos^2\ttheta + 4 \sin^2\ttheta \cos^2\psi),
\]
from which $\cos\psi=0$.  Because of (\ref{range}),
 $\psi=\pi/2$. 
Conversely, when $\psi=\pi/2$, (\ref{eq1})--(\ref{eq3}) are satisfied.
Thus,
\[
 r = 0, \quad \frac{1}{2}\cos^{-1} 0 = \frac{\pi}{4} > \theta_c^{\rm loc}.
\]

(ii)\ \ The case $\sin\psi= 0$.  Then $\cos\psi=\pm 1$ ($\psi=0,\pi$).
In this case both (\ref{eq2}) and (\ref{eq3}) are reduced to
$\sin(\theta \mp \ttheta) = 0$.
Because of (\ref{range}),
$\theta = \pm \ttheta$.
If $\psi=0$ and $\theta = \ttheta$, then $(h,\theta)=(\th,\ttheta)$,
or $x=\tx$.  Hence, it should be
$\psi=\pi$, $\theta = -\ttheta$, and
\[
 r = -1, \quad \frac{1}{2}\cos^{-1} (-1) = \frac{\pi}{2} > \theta_c^{\rm loc}.
\]
We have proved that the critical radius is attained locally.

\subsection{Proof of the recurrences (\ref{forward}) and (\ref{backward})}
\label{subsec:recurrence}


For
$
v = v(\theta)
= 3 \cos^2\theta + 4\sin^2\theta
= 4 - \cos^2\theta = 3 + \sin^2\theta
$,
\begin{eqnarray*}
E_k
&=& \int_{-\pi/2}^{\pi/2} v(\theta)^k d\theta \\
&=& \int_{-\pi/2}^{\pi/2} (4-\cos^2\theta) v^{k-1} d\theta
= 4 E_{k-1} - \int_{-\pi/2}^{\pi/2} \cos\theta^2 v^{k-1} d\theta \\
&=& 4 E_{k-1} - \sin\theta \cos\theta v^{k-1} \bigg|_{-\pi/2}^{\pi/2}
 + \int_{-\pi/2}^{\pi/2} \sin\theta \{ \cos\theta v^{k-1} \}' d\theta \\
&=& 4 E_{k-1} - \int_{-\pi/2}^{\pi/2} \sin^2\theta v^{k-1} d\theta
 + \int_{-\pi/2}^{\pi/2} \sin\theta \cos\theta (k-1) v^{k-2} 2\sin\theta\cos\theta d\theta \\
&=& 4 E_{k-1} - \int_{-\pi/2}^{\pi/2} (v-3) v^{k-1} d\theta
 + 2(k-1) \int_{-\pi/2}^{\pi/2} (v-3)(4-v) v^{k-2} d\theta \\
&=& 4 E_{k-1} - E_k + 3 E_{k-1} - 2(k-1) E_k + 14(k-1) E_{k-1} - 24(k-1) E_{k-2}\\
&=& (-2k+1) E_k + (14k-7) E_{k-1} - 24(k-1) E_{k-2},
\end{eqnarray*}
and hence
\[
 2k E_k = 7(2k-1) E_{k-1} - 24(k-1) E_{k-2}
\]
or
\[
 - 2(k+2) E_{k+2} + 7(2k+3) E_{k+1} = 24(k+1) E_k.
\]



\end{document}